\documentclass[11pt,reqno]{amsart}
\usepackage{amsmath,amsfonts,a4,epsfig}
\usepackage{graphicx}
\usepackage[latin1]{inputenc}
\usepackage[french, english]{babel}
\usepackage{color}
\usepackage{enumerate}
\usepackage{setspace}
\usepackage{fancyhdr}

\def\preuve{\begin{proof}}

\newtheorem{defi}{Definition}[section]

\newtheorem{defini}{Definitions}[section]

\newtheorem{obs}{Observation}[section]
\newtheorem{lemm}{Lemma}[section]
\newtheorem{prop}{Proposition}[section]
\newtheorem{rem}{Remark}[section]

\newtheorem{coro}{Corollary}[section]
\newtheorem{theo}{Theorem}[section]
\newtheorem{exem}{Example}[section]

  \newenvironment{demo}{\noindent {\it Proof:}
      \begin{quotation}\noindent}{\end{quotation}\hfill$\square $}

\begin{document}

\title{On the eigenvalues of weighted directed graphs}
\begin{abstract}
We consider a non self-adjoint Laplacian on a directed graph with non symmetric edge weights. We analyse spectral properties of this Laplacian under a Kirchhoff's assumption. Moreover we establish isoperimetric inequalities in terms of the numerical range to show the lack of essential spectrum of Laplacian on \textit{heavy ends} directed graphs. We introduce a special self-adjoint operator and compare its essential spectrum with that of the non self-adjoint Laplacian considered.
\end{abstract}
\author{ MARWA BALTI}
\address{Universit\'e de Carthage, Facult\'e des Sciences de Bizerte: Math\'ematiques et Applications (UR/13ES47) 7021-Bizerte (Tunisie)\\
Universit\'e de Nantes, Laboratoire de Math\'ematique Jean Lauray, CNRS, Facult\'e des Sciences, BP 92208, 44322 Nantes, (France).
}
\email{ balti-marwa@hotmail.fr}



\maketitle
\tableofcontents

\section*{Introduction}
 This article follows up on the ideas in \cite{I.CH} and \cite{R.P} on monotonicity eigenvalues which are relative
to continuous domains. The purpose of this work is to explore in the case of weighted directed graphs, some familiar
 facts of monotonicity proved on domains of $\mathbb{R}^n$ and on compact Riemannian manifolds. Specially, our main question is:
  "can one study the behavior of the eigenvalues of a special operator under perturbations on finite graphs?"
 First, we consider a finite, directed and connected graph $G$ with non symmetric edge weights. Then, we introduce
   the associated non symmetric Laplacian $\Delta_G$. We introduce the self-adjoint operator $S_G=\Delta_G+\Delta^*_G$, such that it
    is easier to examine its spectrum thanks to selfadjointness.
    We give some spectral properties of $\Delta_G$ and we show that the real part of its eigenvalues  can coincide with the eigenvalues of $\dfrac{1}{2}S_G$. Secondly, we study the monotonicity of eigenvalues relative to vertices or edges of $G$. We prove that the $k^{th}$ eigenvalue $\lambda_k$ of $S_G$ is decreasing only in a class of graphs called here \textit{flower-like-graphs}, and it is monotone
     increasing in the set of edges. These results are inspired by the classical results going back from M. Fiedler
     \cite{M.F} and P. Kurasov, G. Malenova, S. Naboko \cite{Kuras} on the first nonzero eigenvalue of a simple graph.
    We extend these results for the higher eigenvalues $\lambda_k$ of our special operator. In the second part
   of our work, we try to establish and improve the Proposition 2.1 of \cite{B.B} for a Riemannian manifold $M$ which
 gives upper bounds on the higher eigenvalues $\lambda_k(M)$ in terms of Dirichlet eigenvalues on components of a
 partition of $M$.

Let us briefly outline the contents of this article. We shall start with a short section of preliminaries consisting
 on some basic properties of the non symmetric Laplacian $\Delta_G$ on $G$ and the associated Green formula.
In section \ref{monot}, we establish a generalization of some monotonicity eigenvalue results. Furthermore, we study
 the example of flower-like graphs and insist on the interest of the subgraph concept. Hence, we can remark that
 these considerations of graph help to give an upper bound of the eigenvalues of a simple tree. This section includes
     also similar Weyl and Cauchy theorems for the matrices \cite{H.J}.
In section \ref{comp}, we are interested on the study of the eigenvalues of the Dirichlet Laplacian.
We involve a comparison on the eigenvalues to use the decomposition of $G$ into two components $A$ and $B$ and
 we give an upper bound on the eigenvalues of $G$ in terms of Dirichlet eigenvalues of $A$ and $B$.

\section{Preliminaries}
We will review in this section some basic definitions and introduce the notation used in the article.
 They are introduced in \cite{bal} for the infinite graph.
\subsection{Notion of Graphs}
 We call oriented or directed graph, the couple $G=(V,\vec{E})$, where $V$ is a set of vertices,
 and $\vec{E}\subset V\times V$ is a set of directed edges. For two vertices $x,y$ of $V$,
  we denote by $(x,y)$ the edge that connects $x$ to $y$, we also say that $x$ and $y$ are \textit{ neighbors}.
 
For all $x\in V$, we set:
\begin{itemize}
\item $E=\left\{\{x,y\},~(x,y)\in\vec{E}\text{ or } (y,x)\in\vec{E}\right\}$
  \item $V_x^+=\left\{y\in V,~~(x,y)\in\vec{E}\right\}$
  \item $V_x^-=\left\{y\in V,~~(y,x)\in\vec{E}\right\}$
  \item $V_x=V_x^+\cup V_x^-$.

     \end{itemize}
  The valency on $G$ is given by:
  $$v(x)=\#V_x~~\text{ for all } x\in V.$$

We introduce some definitions given in \cite{CN}, \cite{bal}, \cite{R.D}, \cite{N.T}  for the case of symmetric graphs.
\begin{defini}
\begin{itemize}
\item A path between two vertices $x$ and $y$ in $V$ is a finite set of directed edges
$(x_1,y_1);~(x_2,y_2);..;(x_n,y_n),~n\geq 1$ such that
$$x_1=x,~y_n=y \text{ and } x_i=y_{i-1}~~\forall~2\leq i\leq n.$$
\item $G=(V,\vec {E})$ is called connected if two vertices are always related by a path.
\item $G=(V,\vec {E})$ is called strongly connected if there is for all vertices $x,y$ a path from $x$ to $y$
and one from $y$ to $x$.
\begin{exem}
The cycle graph $C_n=\{0,1,..,n-1\}$, with
$$(0,1)\in\vec{E},~(1,2)\in\vec{E},..,(n-2,n-1)\in\vec{E},~(n-1,0)\in\vec{E}$$
is strongly connected.
\end{exem}

\item Define for a finite subset $\Omega$ of ~$V$, the interior, the vertex boundary and the edge boundary
of $\Omega$ respectively by:
$$\overset{\circ}{\Omega}=\big\{y\in \Omega,~ V_y \subset\Omega\big\}$$
$$\partial\Omega=\big\{y\in\Omega^c,~y\in V_x~\text{ for some } x\in \Omega\big\}$$
$$\partial_E \Omega=\big\{(x,y)\in \vec{E}:~ (x\in \Omega,~y\in \Omega^{c} )~~or ~~(x\in \Omega^{c},~y\in \Omega)\big\}.$$
 \end{itemize}
\end{defini}

We remark that the strong connectedness of $G$ assures that:
\begin{equation}\label{cnx}
\forall ~x\in V, ~~\#V_x^+\neq 0~and~\#V_x^-\neq 0.
\end{equation}
In this work we suppose that $G$ is \textbf{finite}, \textbf{connected} and satisfies the  \textbf{Hypothesis} (\textbf{\ref{cnx}}). In this work, we take the following definition.
\begin{defi}\label{an}
\textbf{Directed weighted Graph}: A weighted graph $(G,b)$ is the data of a graph $(V,\vec{E})$ and a weight
$b:V\times V\to \mathbb{R}_+$ satisfying the following conditions:
\begin{itemize}
  \item $b(x,x)=0$ for all $x\in V$, ~~(no loops in  $\vec{E}$)
  \item $b(x,y)>0$ iff $(x,y)\in \vec{E}$
  \item Assumption $(\beta)$.
  \end{itemize}
  \end{defi}
\textbf{Assumption $(\beta)$}: for all $x\in V$,~ $\beta^+(x)=\beta^-(x)$\\
  where $$\beta^+(x)=\sum_{y\in V_x^+ }b(x,y)\text{ and }\beta^-(x)=\sum_{y\in V_x^-}b(y,x).$$
 The weight $\beta_G$ on a vertex $x\in V$ is given by:
$$\beta_G(x)=\beta^+(x)+\beta^-(x)=2\beta^+(x).$$
  \begin{rem}
The Assumption $(\beta)$ is natural, it looks like the Kirchhoff's law in the electrical networks.
\end{rem}
  The weighted graph is \textit{symmetric} if  for all $x,y\in V$, $b(x,y)=b(y,x)$,\\
  as a consequence $(x,y)\in \vec{E}\Rightarrow (y,x)\in\vec{E}$ (the graph is symmetric).\\
In addition, we consider a  weight $m$ on $V$:
$$m:V\to \mathbf{R}^*_+.$$
\subsection{Functional spaces}
Let us introduce the following function spaces associated to the graph $G$:
$$\mathcal{C}_m(V)=\{f:V\to \mathbb{C} \}$$
 endowed with the following inner product:
$$(f,g)_m=\sum_{x\in V}m(x)f(x)\overline{g(x)}.$$
We define its associated norm by:
$$\|f\|_m=\sqrt{(f,f)_m}.$$
A particular case called normalized is for $m = \beta^+$.\\
 For a subset $U$ of $V$, Let
$$\mathcal{C}_m(U)=\{f\in \mathcal{C}_m(V),~ f \text{ with support in U}\}.$$
The weights $m$ and $b$  are called simple if they are constant equal to $1$ on $V$ and $E$ respectively.
 We denote by $G^{s}$ the simple graph (with simple weights).
\section{Laplacian on directed graphs}
For a weighted connected directed graph $(G,b)$, we introduce the combinatorial Laplacians:
\begin{defini}
\begin{itemize}
  \item We define the Laplacian $\Delta_G$ on $\mathcal{C}_m(V)$ by:
$$\Delta_G f(x)=\frac{1}{m(x)}\sum_{y\in V_x^+}b(x,y)\big(f(x)-f(y)\big).$$
  \item In particular, if for all $x\in V$, $\beta^{+}(x)=m(x)$, the Laplacian is said to be the
  normalized Laplacian and is defined on  $\mathcal{C}_{\beta^+}(V)$ by:
 $$\tilde{\Delta}_Gf(x)=\frac{1}{\beta^+(x)}\sum_{y\in V^+_x}b(x,y)\big(f(x)-f(y)\big).$$
  \item For any operator $A$ on $\mathcal{C}_m(V)$, the Dirichlet operator $A^{D}_U$,
  where $U $ is a subset of $ V$, is defined by:
   $$f \text{ is with support in }U,~~ A^{D}_U(f)=A(f)\vert_ U.$$
\end{itemize}
\end{defini}

Thanks to Hypothesis $(\beta)$ the adjoint of $\Delta$ has a simple expression.
\begin{defi}
\textbf{Adjoint of an operator}: The adjoint operator $\Delta^*$ of $\Delta$ is defined by:
 $$\forall \phi,\psi \in \mathcal{C}_m(V), ~( \Delta \psi, \phi) = (\psi,  \Delta^*\phi).$$
 \end{defi}
\begin{prop} Let $f$ be a function of $\mathcal{C}_m(V)$, we have
$$\Delta^*_G f(x)=\frac{1}{m(x)}\sum_{y\in V_x^-}b(y,x)\big(f(x)-f(y)\big).$$
\end{prop}
\begin{demo}
The following calculation for all $f,g \in \mathcal{C}_m(V)$ gives:
\begin{align*}
(\Delta_G f,g)_m=&\sum_{(x,y)\in \vec{E}}b(x,y)\big( f(x)-f(y)\big) \overline{g(x)}\\
=&\sum_{x\in V}\overline{g(x)}f(x)\sum_{y\in V^{+}_x}b(x,y)-\sum_{(y,x)\in \vec{E}}b(y,x) \overline{g(y)}f(x)\\
=&\sum_{x\in V}\overline{g(x)}f(x)\sum_{y\in V^{-}_x}b(y,x)-\sum_{(y,x)\in \vec{E}}b(y,x) \overline{g(y)}f(x)\\
=&\sum_{x\in V} f(x)\sum_{y\in V_x^-}b(y,x)\big( \overline{g(x)-g(y)\big)}\\
=&(f,\Delta^*_Gg)_m.
\end{align*}
\end{demo}

 The Green's formula is one of the main tools when we are working with the symmetric Laplace operator.
 In the following we establish it for the non symmetric Laplacian.
\begin{lemm} \textbf{Green's Formula}.
Let $f$ and $g$ be two functions of $\mathcal{C}_m(V)$. Then
$$(\Delta_G f,g)_m+(\Delta^*_Gf,g)_m=\sum_{(x,y)\in \vec{E}}b(x,y)\big( f(x)-f(y)\big) \big( \overline{g(x)-g(y)}\big).$$
\end{lemm}
\begin{demo}
The proof is a simple calculation:
\begin{align*}
(\Delta_G f,g)_m+(\Delta^*_Gf,g)_m=&\sum_{(x,y)\in \vec{E}} b(x,y)\big(f(x)-f(y)\big)\overline{g(x)}\\
+&\sum_{(y,x)\in \vec{E}} b(y,x)\big(f(x)-f(y)\big)\overline{g(x)}\\
=&\sum_{(x,y)\in \vec{E}} b(x,y)\Big(f(x)\overline{g(x)}+f(x)\overline{g(x)}-f(y)\overline{g(x)} -f(x)\overline{g(y)}\Big)\\
=&\sum_{(x,y)\in \vec{E}} b(x,y)\big( f(x)-f(y)\big) \big(\overline{g(x)-g(y)}\big).
\end{align*}
\end{demo}
\begin{defi}{\textbf{Special Laplacian.}}
We define a special Laplacian $S_G$ as the sum of the two non self-adjoint Laplacians $\Delta_G$ and $\Delta^{*}_G$, given by:
\begin{align*}
S_Gf(x)&=(\Delta_G+\Delta^{*}_G)f(x)\\
&=\frac{1}{m(x)}\sum_{y\in V_x^{+}\cup V_x^{-}}\big(b(x,y)+b(y,x)\big)\big(f(x)-f(y)\big)\\
&=\frac{1}{m(x)}\sum_{y\in V_x}a(x,y)\big(f(x)-f(y)\big)
\end{align*}
where $a(x,y)=b(x,y)+b(y,x)$ for any $x,y\in V$.
\end{defi}
\begin{rem}
\begin{enumerate}
 \item $S_G$ is a symmetric operator on $\mathcal{C}_m(V)$, because $(\Delta^{*}_G)^*=\Delta_G$.
\item $S_G$ is a positive operator: for all $f\in \mathcal{C}_m(V)$, \begin{align*}
                                                                (S_Gf,f)= &(\Delta_Gf,f)+(\Delta_G^*f,f) \\
                                                                = & \sum_{(x,y)\in \vec{E}} b(x,y)\big| f(x)-f(y)\big|^2 \\
                                                                \geq & 0.
                                                              \end{align*}
 \end{enumerate}
\end{rem}
In our discussion on the study of eigenvalues of a self-adjoint operator, it is natural to introduce
 the different characterization by variational principles \cite{Has}.
\subsection{Variational principles and Properties}
Let $A$ be a bounded from below self-adjoint operator. The eigenvalues of $A$ can be characterized by three fundamental variational principles:
 the Rayleigh's principle, the Poincar\'e-Ritz
max-min principle and the Courant-Fischer-Weyl principle applied to the Rayleigh
quotients $\mathcal{R}(f)=\dfrac{(Af,f)}{(f,f)}$, $f\neq 0$.
\\
Let us arrange the eigenvalues of $A$ as
$$\lambda_1\leq \lambda_2\leq...\leq\lambda_n$$
counted according to their multiplicities. \\
In this case we have :
\begin{enumerate}
  \item The \textit{Rayleigh's principle} states:
  \begin{equation}
\lambda_k=\min_{f\neq0,~(f,f_i)=0 \atop i=1,...,k-1}\mathcal{R}(f)
\label{rai}
\end{equation}
where $f_i$ are eigenvectors corresponding to the eigenvalues $\lambda_i$ and the minimum is reached
at the eigenvector $f_k$.
  \item The \textit{Poincar\'e-Ritz principle} establishes:
\begin{equation}
\lambda_k=\min_{dim \Omega=k}~\max_{f\in \Omega,~f\neq 0}\mathcal{R}(f).
\label{poin}
\end{equation}
  \item The \textit{Courant-Fischer-Weyl principle} is given in the form:
\begin{equation}
\lambda_k=\max_{dim \Omega=k-1}~\min_{f\perp \Omega\atop~f\neq 0}\mathcal{R}(f).
\label{fish}
\end{equation}
\end{enumerate}

The result below establishes a link between the eigenvalues of $S_G$ and $\Delta_G$.
 We assume that the eigenvalues of $\Delta_G$ are ordered as follows respectively:
$$\mathcal{R}e(\lambda_1(\Delta_G))\leq\mathcal{R}e(\lambda_2(\Delta_G))..\leq\mathcal{R}e(\lambda_n(\Delta_G)).$$

\begin{lemm}\label{delta}
  $$2\mathcal{R}e(\lambda_n(\Delta_G))\leq\lambda_n(S_G)$$
\end{lemm}
\begin{demo}
Let $f$ be an eigenfunction associated to $\lambda_n(\Delta_G)$, we have
\begin{align*}
  \lambda_n(S_G)\geq&\dfrac{(S_G f,f)_m}{(f,f)_m}\\
  \geq & \dfrac{(\Delta_G f,f)_m+\overline{(\Delta_G f,f)}_m}{(f,f)_m} \\
  \geq & 2\mathcal{R}e(\lambda_n(\Delta_G)).
\end{align*}
\end{demo}
\begin{rem}
In a particular case, the previous inequalities are strict. Let us consider the following example, where $m\equiv 1$,
\begin{figure}[ht]
\begin{center}
\includegraphics*[height=2.5cm,width=3cm]{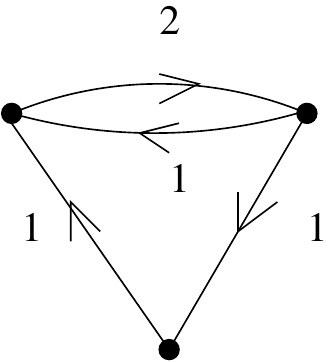}
\end{center}
 \caption{Graph with real spectrum}
\end{figure}
we have  $\sigma(S_G)=\{0,3,7\}$ and $\sigma(\Delta_G)=\{0,2,3\}$.
\end{rem}

We introduce in the following a particular case of graphs whose $\sigma(S_G)=2\mathcal{R}e\big(\sigma(\Delta_G)\big)$.
\begin{exem}
Let us consider the simple cycle graph $C_3$, see the Figure \ref{cyc}, we have ,
$$\sigma(\Delta_{C_3})=\left\{0,~\frac{3}{2}+i\frac{\sqrt{3}}{2},~\frac{3}{2}-i\frac{\sqrt{3}}{2}\right\}\text{ and } \sigma(S_{C_3})=\left\{0,3,3\right\}.$$
\begin{figure}[ht]
\begin{center}
\includegraphics*[height=2cm,width=3cm]{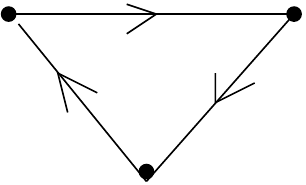}
\end{center}
 \caption{Cycle graph}
 \label{cyc}
\end{figure}
\end{exem}

In the following proposition we determine the spectrum of the non symmetric Laplacian $\Delta_{C_n}$.   We follow the same approach as Grigoryan for the symmetric Laplacian \cite{A.G} page 49.
\begin{prop}
The eigenvalues of $\Delta_{C_n}=I-P$  are as follows:
\begin{enumerate}
\item If $n$ is odd then the eigenvalues are $\lambda=0$ (simple) and $\lambda=1-e^{\pm i\frac{2l\pi }{n}}$ for all $l=1,..,\frac{n-1}{2}$ (simple).
\item If $n$ is even then the eigenvalues are $\lambda=0,~2$ (simple) and $\lambda=1-e^{\pm i\frac{2l\pi}{n}}$ for all $l=1,..,\frac{n}{2}-1$ (simple).
\end{enumerate}
\end{prop}
\begin{demo}
To compute the eigenvalues of $\Delta_{C_n}=I-P$, it is sufficient to determine the spectrum of $P$. Let $\alpha$ be an eigenvalue of the operator  $Pf(k+1)= f(k)$, for $k=0,..,n-1$, which leads to $f(k)=\alpha^k f(0)$ but $f(n)=f(0)$ thus $\alpha^n=1$.  As $f$ is an eigenfunction $(f(0)\neq 0)$, then $\alpha=e^{\pm i\theta}$. As $f$ is n-periodic provided $n\theta$
is a multiple of $2\pi$, hence, $$\theta=\frac{2l\pi}{n},$$
where $l$ is an integer of $(0,\frac{n}{2})$.
\end{demo}

Observe that an interesting corollary concerning the spectra of $S_{C_n}$ and $\Delta_{C_n}$.
\begin{coro}
$$\sigma(S_{C_n})=2\mathcal{R}e\big(
\sigma(\Delta_{C_n})\big).$$
\end{coro}
\begin{demo}
We refer the Lemma 2.7 of \cite{A.G}, we remark that the eigenvalues of the operator $\dfrac{1}{2}(P+P^*)f(k)=\dfrac{1}{2}\big(f(k+1)+f(k-1)\big)$ coincide with the real part of the eigenvalues of $P$.
\end{demo}

Using the Green's formula, we establish some properties of the spectrum on any graph $G$ .
\begin{prop}
\begin{enumerate}
\item $0$ is a simple eigenvalue of $\tilde{S}_G$ and $S_G$.
\item All the eigenvalues of $\tilde{S}_G$ are contained in $[0,4]$.
\item  The real part of the eigenvalues of $\tilde{\Delta}_G$ are also contained in $[0,2]$.
\end{enumerate}
\end{prop}
\begin{demo}
\begin{enumerate}
\item As in the case of undirected graph \cite{A.G}, we have for all $f\in\mathcal{C}_{\beta^+}(V)$,
$$(\tilde{S}_Gf,f)_{\beta^+}=\sum_{(x,y)\in \vec{E}} b(x,y)|f(x)-f(y)|^2.$$
Clearly, the constant function is an eigenfunction of $0$. Assume now that f is an eigenfunction
 of the eigenvalue $0$. By the connectedness of $G$, $f$ is constant, which will imply that $0$ is a simple eigenvalue.
 It is similar for $\tilde{S}_G$.
\item It is sufficient to prove that $\tilde{S}_G$ is bounded by $4$ because $\tilde{S}_G$ is non negative by the Green's formula. In fact, for all $f\in \mathcal{C}_{\beta^+}(V)$
 and thanks to Assumption $(\beta)$, we obtain
\begin{align*}
(\tilde{S}_Gf,f)_{\beta^+}=&\sum_{(x,y)\in \vec{E}} b(x,y)|f(x)-f(y)|^2\\
\leq &2\sum_{(x,y)\in \vec{E}} b(x,y)\big||f(x)|^2+|f(y)|^2\big|\\
\leq &2\sum_{(x,y)\in \vec{E}} b(x,y)|f(x)|^2+2\sum_{(x,y)\in \vec{E}} b(x,y)|f(y)|^2\\
\leq &2\sum_{x\in V}\sum_{y\in V_x^+} b(x,y)|f(x)|^2+2\sum_{y\in V}\sum_{x\in V_x^-} b(x,y)|f(y)|^2\\
\leq &2\sum_{x\in V}|f(x)|^2\beta^+(x)+2\sum_{y\in V}|f(y)|^2\beta^-(y)\\
\leq & 4(f,f)_{\beta^+}.
\end{align*}
\item We deduce directly our inclusion thanks to the Lemma \ref{delta}.
\end{enumerate}
\end{demo}
\section{Domain monotonicity of eigenvalues}\label{monot}
 The purpose of this part is to give an overview of some  results concerning the monotonicity with regard to
 the domain, of eigenvalues of $S_G$, the special self-adjoint Laplacian associated to directed graphs
 with non symmetric edge weights. We could be concerned with the related question:\\
Does a given eigenvalue increases or decreases under a given perturbation of $G$?
\subsection{Definitions on $G$}
  Before discussing the study of variation of eigenvalues, let us recall some basic definitions:
   let $G=(V,\vec{E})$ be a graph,
\begin{itemize}
  \item The graph $G_1=(V, \vec{E}_1)$ is called a partial graph of $G$, if $\vec{E}_1$ is included in $\vec{E}$.
  \item A graph $H=(V_H,\vec{E}_H)$ is called a subgraph of $G=(V_G,\vec{E}_G)$ if $V_H\subset V_G$
and $\vec{E}_H=\big\{(x,y) ;~x,y \in V_H~~\big\}\cap \vec{E}_G $.
  \item A graph $(V_U,\vec{E}_U)$ is called a part of a graph $G=(V_G,\vec{E}_G)$ if $V_U\subset V_G$ and
  $\vec{E}_U=\big\{(x,y),~~x,y\in V_U\big\}\subset \vec{E}_G.$
\end{itemize}
\begin{rem}
A subgraph is a part of $G$ but the converse is not true, for example let us give the following undirected graphs,
 see the Figure \ref{exp}.
\begin{figure}[ht]
\begin{center}
\includegraphics*[height=2cm,width=9cm]{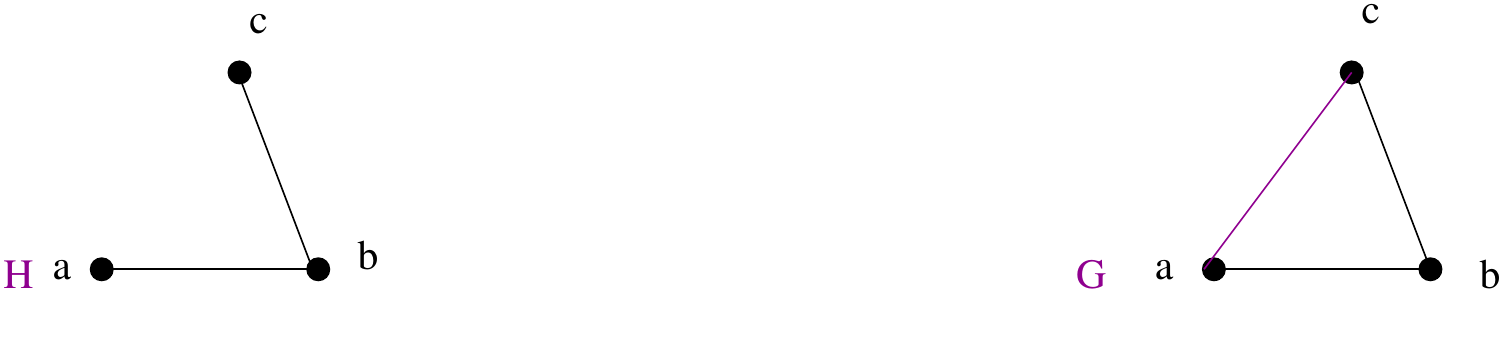}
\end{center}
 \caption{$H$ part of $G$}
 \label{exp}
\end{figure}
\\
 $H$ is a part of $G$, but not a subgraph.
 \end{rem}
 \begin{rem}
A confusion between a subgraph and a part of a graph can create a false interpretation on the
 monotonicity of eigenvalue.
\end{rem}
  \subsection{Monotonicity relative to vertices } We study the monotonicity of eigenvalues under the variation of the set of vertices.\\

By the Courant-Fischer-Weyl principle we establish the following statement.
\begin{theo}\label{thm1}
Let $H=(V_H,\vec{E}_H)$ be a connected subgraph of a graph  $G=(V_G,\vec{E}_G)$, $\sharp V_G=n$, then for any
 $1 \leq k\leq \sharp V_H=r$ :
$$\lambda_k(S_H) \leq \lambda_{n-r+k}(S_G).$$
\end{theo}
\begin{demo}
Let $f_1,~f_2,...,f_{k}$ be $k$ eigenfunctions associated to $\lambda_1(S_H),~\lambda_2(S_H),...,\lambda_{k}(S_H)$;
and $F=\{f_1,~f_2,...,f_{k},\delta_1,..,\delta_{n-r}\}$, $\delta_1,..,\delta_{n-r}$ are the Dirac measures on $G$ relative to the vertices in $V_G\setminus V_H$.
It is clear that $\dim F=k+n-r$.
 Then using (\ref{fish}), we obtain ;

  $$ \lambda_{n-r+k+1}(S_G) \geq \min_{\varphi\in F^\perp\setminus\{0\}}\dfrac{(S_G \varphi,\varphi)_m
   }{(\varphi,\varphi)_m }$$
   hence $\exists\varphi_k\in F^\perp$ with support in $H$ such that
   \begin{align*}
    \lambda_{n-r+k+1}(S_G)( \varphi_k,\varphi_k)_m &\geq(S_G \varphi_k,\varphi_k)_m\\
    &=(S_H \varphi_k,\varphi_k)_m\\
      &\geq \lambda_{k+1}(S_H)(\varphi_k,\varphi_k)_m.
 \end{align*}
\end{demo}

For studying the behavior of eigenvalues relative to perturbations, we propose a special construction of graphs.
\begin{defi}
Let $G=(V_G,\vec{E}_G)$ be a weighted graph. $G$ is called $H$-flower-like with respect to the
subgraph $H=(V_H,\vec{E}_H)$ of $G$ if there exists $(H_i)_{i\in I}$ a family
 of subgraphs of $G$ such that:
\begin{enumerate}
 \item $V_G=V_{H}\cup (\displaystyle{\uplus_i  V_{\overset{\circ}{{H_i}}}})$
  \item $\forall ~i,j\in I,~~~x\in V_{\overset{\circ}{{H_i}}},~~y\in V_{\overset{\circ}{{H_j}}}, ~i\neq j~\Rightarrow~x\notin V_y$.
  \item $\forall~i\in I,~\exists ~x_i\in V_G,~V_H\cap V_{H_i}=\{x_i\}$.
\end{enumerate}
\end{defi}
\begin{figure}[ht]
\begin{center}
\includegraphics*[height=6cm,width=7cm]{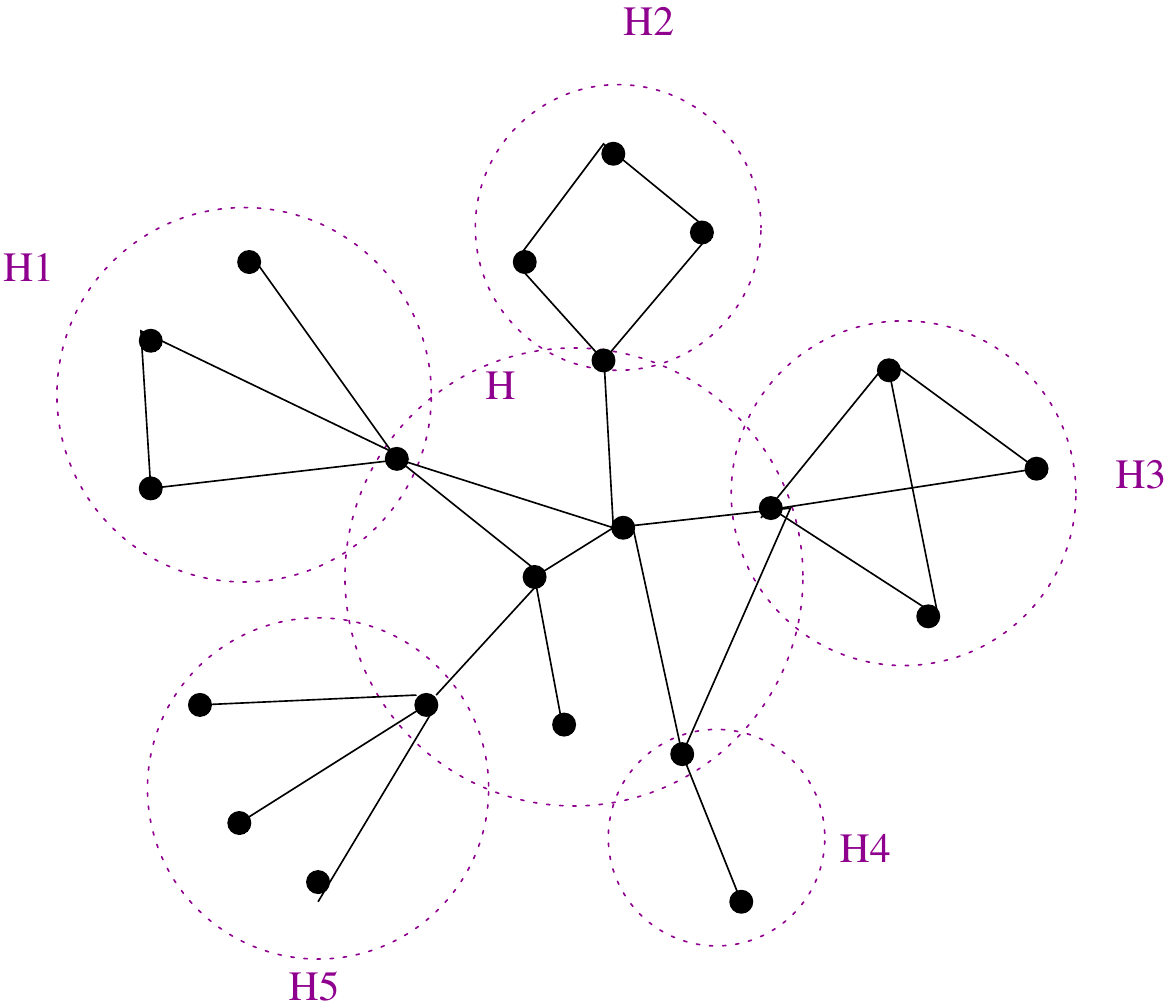}
\end{center}
 \caption{H-flower-like graph}
\end{figure}
These special graphs are used to create a rule of monotonicity of $\lambda_k$ under a given graph perturbation.
 \begin{theo} \label{AA}
 Let $G=(V_G,\vec{E}_G)$ be a $H$-flower-like graph, we have then for any $1\leq k\leq \sharp V_H=r$,
 \begin{equation}\label{Hfl}
 \lambda_k(S_H) \geq \lambda_k(S_G).
 \end{equation}
\end{theo}
\begin{demo}
We use the variational principle (\ref{rai}). Let be $f_1,~f_2,...,f_{k-1}$ eigenfunctions of $G$ in $\mathcal{C}_m(V)$
 associated to $\lambda_1(S_G),~\lambda_2(S_G),...,\lambda_{k-1}(S_G)$,
and $g_j$ be the eigenfunction associated to $\lambda_j(S_H)$, for $j=1,..,k$. We define a function $\phi_j$ on $V_G$ by:
\begin{equation*}
\phi_j (x)=
  \left\{
        \begin{aligned}
        & g_j (x)~~~~~~~~~~~~~~\text{ if } x\in V_H\\
        & g_j (x_1)~~~~~~~~~~~~\text{ if } x\in V_{H_1}\setminus\{x_1\}\\
        & . \\
        & . \\
        & g_j (x_r)~~~~~~~~~~~~\text{ if } x\in V_{H_r}\setminus \{x_r\}
      \end{aligned}
    \right.
\end{equation*}
where $\{x_s\}=V_H\cap V_{H_s}, ~\forall s\in I=\{1,..,r\}$\\
Let $F=Span\{\phi_j\}_{1\leq j \leq k}$, $\dim F=k$ so, there exist $k$ reals $\alpha_1,\alpha_2,...,\alpha_k$
not all equal to zero, satisfying:
$$\Big(\sum_{j=1}^{k}\alpha_j \phi_j~,~f_l\Big)_m=0~~~~~\forall~~ 1 \leq l \leq k-1.$$
Therefore the function $\psi(x)=\displaystyle{\sum_{j=1}^{k}}\alpha_j \phi_j(x)\neq 0 $ is orthogonal
to $f_l$ for all $1 \leq l \leq k-1$, we define for $x\in H$, $\psi'(x)=\displaystyle{\sum_{j=1}^{k}}\alpha_j g_j(x)$.
Then we get :
\begin{align*}
   \lambda_k(S_G)(\psi,\psi)_m & \leq (S_G \psi,\psi)_m \\
   &=(S_H \psi',\psi')_m ~~\text{ because $\psi$ is constant on each $H_i$ }\\
      &\leq\lambda_k(S_H)(\psi',\psi')_m
\end{align*}
as
\begin{align*}
  (\psi,\psi)_m~=~ & \sum_{x\in V_G}m(x)|\psi(x)|^{2}\\
  \geq~& \sum_{x\in V_H}m(x) |\psi'(x)|^{2}\\
  =~&(\psi',\psi')_m
\end{align*}
then
$$\lambda_k(S_G)(\psi',\psi')_m  ~~\leq ~~\lambda_k(S_H)(\psi',\psi')_m.$$
\end{demo}

The above results have several important consequences, for instance on a tree seen as a flower-like graph.
\begin{coro}
Let $G^{s}$ be a simple symmetric tree and $2q=\displaystyle{\max_{x\in V}v(x)}$ then the eigenvalues
of $G^{s}$ satisfy for all $k=1,..,q$:
$$\lambda_k(S_{G^{s}})\leq 2$$
and
$$\lambda_{q+1}(S_{G^{s}})\leq 2(q+1).$$
\end{coro}
\begin{demo}
 Clearly there exists a symmetric star graph $S_q$ with $q+1$ vertices seen as a subgraph of $G^{s}$. It can be
 considered as a $S_q$-flower-like graph. Therefore $G^{s}$ satisfies the assumptions of the Proposition \ref{AA}.
  Hence we get the result because the spectrum of $S_q$ is : $0,2,...,2,2(q+1)$, see \cite{T.M}.
\end{demo}

In the following, we will show more general results : instead of adding only one vertex and one edge,
we would also add a graph.
\begin{coro}
 Let $G$ be a graph with $n$ vertices, let $G_1$ connected to G by a single edge.
Then for $k=1,..,n$;
$$\lambda_k(S_G) \geq \lambda_k(S_{G_1}).$$
\end{coro}
\begin{rem}
The previous corollary is an immediate consequence of the Proposition \ref{AA}. In addition this is an
interesting generalization of Proposition 2 in \cite{Kuras} which shows that $\lambda_2(S_G) \geq \lambda_2(S_{G_1})$
 where $G_1$ is obtained from $G$ by adding one edge between one vertex of $G$ and one new vertex.
\end{rem}
\subsection{Monotonicity relative to edges}
We apply the Weyl and interlacing Theorems for matrices to study how the spectrum of the special
operator $S_G$ of a directed graph changes under adding an edge or a set of edges.\\

Now we have several opportunities to refer to the following basic observation about subspace intersections
 (see \cite{H.J} page 235).
\begin{lemm}\label{lemm}
Let $W$ be a finite dimensional vector space and let $S_1$, $S_2$,..,$S_k$ be subspaces of $W$, if
$$\delta=\dim S_1+...+\dim S_k -(k-1) \dim W \geq 1 $$
 then~~$\dim (S_1\cap ..\cap S_k)\geq \delta \text{ and hence  } S_1\cap ..\cap S_k \text{ contains } \delta$
 linearly  independent vectors, in particular , it contains a nonzero vector.
\end{lemm}
\begin{obs}
Let $S_G$ be the special self-adjoint operator with eigenvalues $\lambda_1(S_G) \leq \lambda_2(S_G)\leq...
\leq\lambda_n(S_G)$. Then the ordered eigenvalues of $-S$  are
$\lambda_1(-S_G) \leq \lambda_2(-S_G)\leq...\leq\lambda_n(-S_G)$, that is, $\lambda_k(-S_G)=-\lambda_{n-k+1}(S_G),~k = 1, . . . , n.$
\end{obs}

  We show that the $k^{th}$ eigenvalue $\lambda_k(S_G)$ is monotonously increasing functions of the set of edges.
The following results are the generalizations of the Proposition 1 and the Proposition 2 in \cite{Kuras}.
\begin{prop} \label{Aaa}
Let $G=(V,\vec{E})$ be a connected finite weighted graph $(\sharp V=n)$, consider the partial
 graphs $G_1=(V,\vec{E}_1)$ and $G_2=(V,\vec{E}_2)$ where $E=E_1\sqcup E_2$ (disjoint union).
 Then for all $k=1,...,n$ and $r,s=1,2$, $r\neq s$:
\begin{equation}\label{equ}
\lambda_k(S_G) \leq \lambda_{k+j}(S_{G_s})+ \lambda_{n-j}(S_{G_r})~~~~j=0,...,n-k
\end{equation}
and
$$\lambda_k(S_G) \geq \lambda_{k-j+1}(S_{G_s})+ \lambda_{j}(S_{G_r})~~~~j=1,...,k.$$
\end{prop}
\begin{demo}
Let $g_k$,~ $h_k$, and $f_k$ be the eigenfunctions associated to $\lambda_k(S_{G_1})$, $\lambda_k(S_{G_2})$
 and $\lambda_k(S_G)$ respectively for $k=1,..,n$. \\
For $j=0,...,n-k$, we define $S_1=Span\{g_1,...,g_{k+j}\}$, $S_2=Span\{h_1,...,h_{n-j}\}$ and $S_3=Span\{f_k,...,f_n\}$
by the Lemma (\ref{lemm}) there exists a non zero function $\psi$ in $S_1\cap S_2 \cap S_3$
and so we will have
\begin{align}\label{pos}
\lambda_k(S_G)(\psi,\psi)_m &\leq (S_{G} \psi,\psi)_m\\
&=\sum_{(x,y)\in \vec{E}}b(x,y)\big|\psi(x)-\psi(y)\big|^{2}\notag\\
&= (S_{G_s} \psi,\psi)_{m}+(S_{G_r} \psi,\psi)_{m}\notag\\
&\leq \lambda_{k+j}(S_{G_s})(\psi,\psi)_{m}+ \lambda_{n-j}(S_{G_r})(\psi,\psi)_{m}\notag\\
&\leq\big( \lambda_{k+j}(S_{G_s})+ \lambda_{n-j}(S_{G_r})\big)(\psi,\psi)_m\notag.
\end{align}
In the following we apply the equality (\ref{equ}) to the operator $-S_{G} $ because the inequalities
 (\ref{pos}) does not depend on positivity of $b$, and by:
 $$\lambda_{k}(-S_G)=-\lambda_{n-k+1}(S_G).$$
  We obtain by re-indexing:
  $$\lambda_k(S_G) \geq \lambda_{k-j+1}(S_{G_s})+ \lambda_{j}(S_{G_r})~~~~j=1,...,k.$$
\end{demo}

We can easily deduce,
\begin{coro} \label{Aa}
Let $G=(V,\vec{E})$ be a connected finite graph, $G_1=(V,\vec{E}_1)$ and $G_2=(V,\vec{E}_2)$ two partial
graphs of $G$ where $\sharp V=n \text{ and }E=E_1\sqcup E_2$, then for all $k=1,...,n$ and
$r,s=1,2$, $r\neq s$:
$$\lambda_{k}(S_{G_s}) \leq  \lambda_k(S_G) \leq \lambda_{k}(S_{G_s})+ \lambda_{n}(S_{G_r}).$$
\end{coro}
\begin{demo}
By applying the proposition \ref{Aaa} to $j=0$ and $j=1$ respectively, we obtain the result because
$\lambda_{1}(S_{G_r})=0.$
\end{demo}

In other words, adding a subset of edges to $\vec{E}$ while keeping the same set of vertices always induces
an increasing of the $k^{th}$ eigenvalue or keeps it unchanged.
\begin{coro}\label{cor}
Let $G=(V,\vec{E})$ be a connected finite graph with $n$ vertices, and $G_1$ a graph obtained by adding a set
of edges to $G$ then for all $k=1,...,n$:
$$\lambda_{k}(S_G) \leq  \lambda_k(S_{G_1}).$$
\end{coro}

\section{Comparison eigenvalues of Dirichlet Laplacian on graphs}\label{comp}
In this section, we present some results about the spectrum comparison between the Laplacian and the Dirichlet Laplacian.
 The purpose of this part is to find the relation between the usual vertex weight on a subgraph $H$ of $G$ and its
  boundary weight to compare eigenvalues.\\
This is done by establishing a clear and explicit link between the eigenvalues and the Dirichlet eigenvalues on $G$.

In the following proposition, we treat the Dirichlet Laplacian case :
$$\forall x\in U,~~S^D_U f(x)=\frac{1}{m(x)}\sum_{y\in V_x\atop y\in V_G}a(x,y)\big(f(x)-f(y)\big).$$
By the same techniques used in the Lemma \ref{delta}, we can show the following Lemma.
\begin{lemm}
  $$\lambda_1(S^D_U)\leq 2\mathcal{R}e(\lambda_1(\Delta^D_U))$$
and
  $$\lambda_n(S^D_U)\geq 2\mathcal{R}e(\lambda_n(\Delta^D_U)).$$
\end{lemm}
\begin{demo}
Let $f$ and $g$ be eigenfunctions associated to $\lambda_1(\Delta^D_U)$ and $\lambda_n(\Delta^D_U)$ respectively.
By the variational principle of $\lambda_1(S^D_U)$ and $\lambda_n(S^D_U)$, we have
$$\lambda_1(S^D_U)\leq\dfrac{(S^D_U f,f)_m}{(f,f)_m}=\dfrac{(\Delta^D_U f,f)_m}{(f,f)_m}+
\dfrac{\overline{(\Delta^D_U f,f)}_m}{(f,f)_m}=\lambda_1(\Delta^D_U)+\overline{\lambda_1(\Delta^D_U)}$$
and
$$\lambda_n(S^D_U)\geq\dfrac{(S^D_U g,g)_m}{(g,g)_m}=\dfrac{(\Delta^D_U g,g)_m}{(g,g)_m}
+\dfrac{\overline{(\Delta^D_U g,g)}_m}{(g,g)_m}=\lambda_n(\Delta^D_U)+\overline{\lambda_n(\Delta^D_U)}$$
\end{demo}

In the same spirit as the Cauchy interlacing theorem concerning hermitian bordered matrices (see [9]
 theorem 4.3.28 for a generalized statement) one can prove the following.

\begin{prop}\label{thm2}
Consider a connected subgraph $H=(V_H,\vec{E}_H)$ of $G=(V_G,\vec{E}_G)$, $(\#V_G=n, \#V_H=r)$,
then the eigenvalues on $G$ of the normalized Laplacian $S_G$ satisfies:
$$\lambda_k(S^D_H) \leq \lambda_{k+n-r}(S_G).$$
\end{prop}
\begin{demo}
Let $h_1,..,h_r$ and $f_1,..,f_n$ be the eigenfunctions associated to $S^{D}_H$ and $S_G$ respectively,
 define the function $g_k$ for $k=1,..,r$ by:
\begin{equation*}
g_k=~\left\{
\begin{aligned}
 &h_k~~~~on~~V_H  \\
 &0~~~~otherwise
  \end{aligned}
\right.
\end{equation*}
Let $1\leq k\leq r$ and fix $S_1=Span\{g_k,...,g_{r}\}$ and $S_2=Span\{f_1,...,f_{k+n-r}\}$
by the Lemma \ref{lemm} there exists a function $\psi$ in $S_1\cap S_2$. Since  $\psi \in S_1$, it has the form
\begin{equation*}
\psi=~\left\{
\begin{aligned}
 &g~~~~on~~V_H  \\
 &0~~~~otherwise
  \end{aligned}
\right.
\end{equation*}
for some $g\in Span\{g_k,...,g_{r}\}$. Observe that:
\begin{align*}
\lambda_k(S^D_H)( g,g)_{m} &\leq(S^{D}_{H} g,g)_{m}\\
&\leq (S_G \psi,\psi)_{m}\notag\\
\end{align*}
since $( g,g)_{m}=(\psi,\psi)_{m}$, we get
$$\lambda_k(S^D_H)\leq\dfrac{(S_G\psi,\psi)_{m}}{(\psi,\psi)_{m}}$$
then,
$$\lambda_k(S^D_H) \leq \lambda_{k+n-r}(S_G).$$
\end{demo}

We deduce easily from  the Proposition \ref{thm2} an estimation of the eigenvalues of $S_G$ and $\Delta_G$,
thanks to Theorem \ref{thm1}, as follows:
\begin{coro}
Consider a connected subgraph $H=(V_H,\vec{E}_H)$ of $G=(V_G,\vec{E}_G)$, $(\#V_G=n, \#V_H=r)$,
then the eigenvalues on $G$ satisfies:
$$ \lambda_{k+n-r}(S_G)\geq \max\big(\lambda_k(S_H),\lambda_k(S^D_H)\big).$$
\end{coro}
\begin{coro}
Consider a connected subgraph $H=(V_H,\vec{E}_H)$ of $G=(V_G,\vec{E}_G)$, $(\#V_G=n, \#V_H=r)$,
then:
$$ \lambda_{n}(S_G)\geq \max\Big(2\mathcal{R}e\big((\lambda_r(\Delta_H)\big),2\mathcal{R}e\big(\lambda_r(\Delta^D_H)\big)\Big).$$
\end{coro}
\begin{coro}
Consider a connected subgraph $H=(V_H,\vec{E}_H)$ of the cycle graph $C_n$, $(\#V_H=r)$,
then:
$$\mathcal{R}e\big((\lambda_n(\Delta_{C_n})\geq \max\Big(\mathcal{R}e\big(
(\lambda_r(\Delta_{H})\big),\mathcal{R}e\big(\lambda_r(\Delta^D_{H})\big)\Big).$$
\end{coro}
 In the following Proposition we prove how to give an upper bound for $\lambda_{k+l}(S_G)$ in
terms of a Rayleigh quotient. We give a discrete version of the Proposition 2.1 \cite{B.B} by applying
the Poincar\'e min-max principle. The methods we use follow closely the arguments given in B. Benson \cite{B.B}
in the case of the Laplacian of Riemannian manifolds.  Next we provide an upper bound of the eigenvalues of $G$
according to Dirichlet eigenvalues on such repartition of $G$ as in the Figure \ref{partition}.
\begin{prop}
Let $G=(V,\vec{E})$ a finite connected graph, $U=(V_U,\vec{E}_U)$ a part of $G$ and $A=(V_A,\vec{E}_A)$,
$B=(V_B,\vec{E}_B)$ two subgraphs satisfying the following conditions:
\begin{enumerate}
  \item $V=V_A\sqcup V_B=V_{\overset{\circ}{A}}\sqcup V_{\overset{\circ}{B}}\sqcup V_U$ ~~~(\text{ disjoint union})
  \item $\vec{E}=\vec{E}_A\sqcup\vec{E}_B\sqcup \vec{E}_U$
  \item $\partial_{\vec{E}} A=\partial_{\vec{E}} B=\vec{E}_U$.
\end{enumerate}
\begin{figure}[ht]
\begin{center}
\includegraphics*[height=5cm,width=7cm]{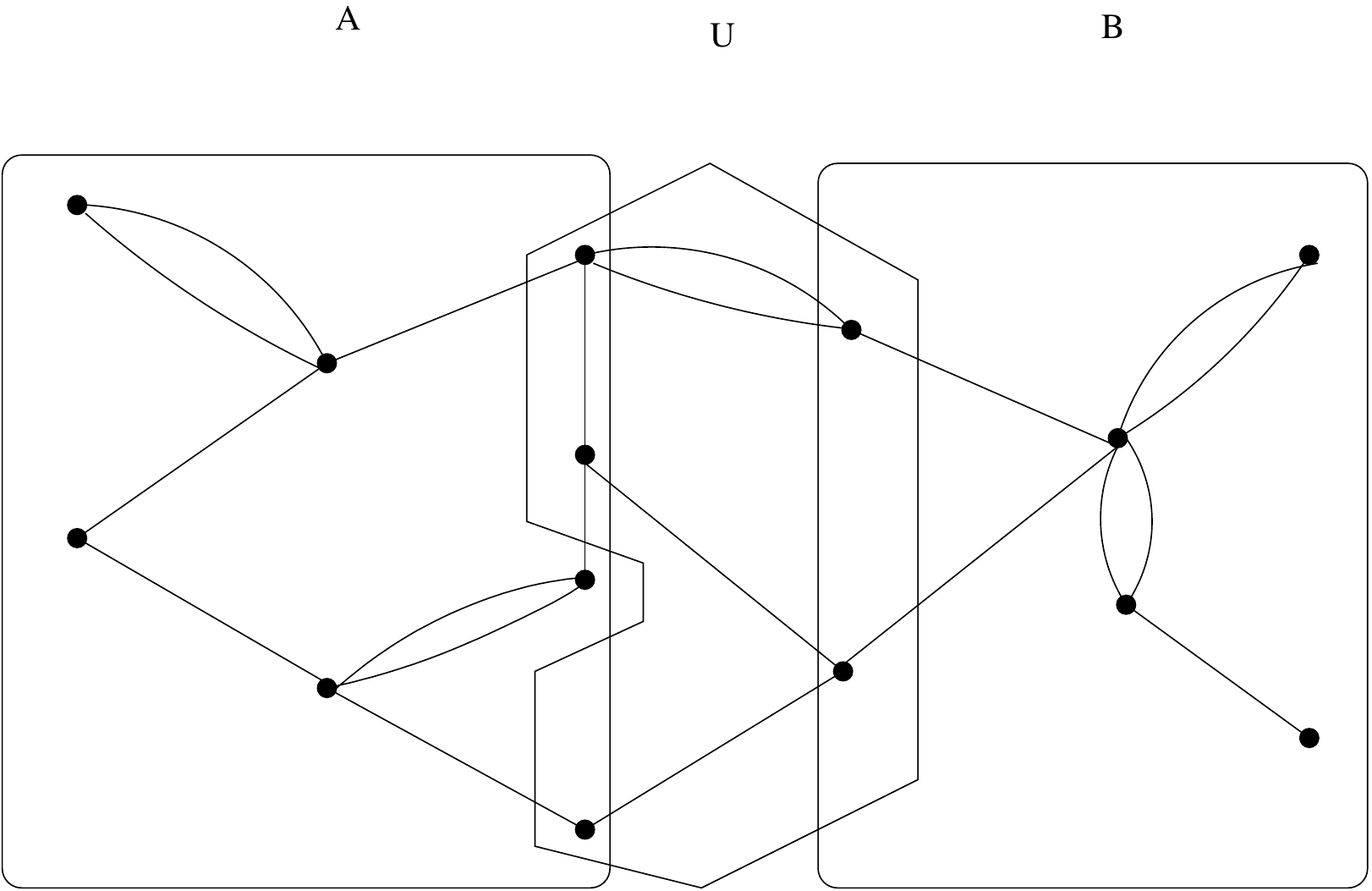}
\end{center}
 \caption{Repartition of $G$}
 \label{partition}
\end{figure}
Then we get for $1 \leq k,l \leq \min(\sharp V_{\overset{\circ}{A}},\sharp V_{\overset{\circ}{B}})$:
$$\lambda_{k+l}(S_G)\leq \max \big(\lambda_{k}(S^D_{\overset{\circ}{A}}),\lambda_l(S^D_{\overset{\circ}{B}})\big).$$
\end{prop}
 
\begin{demo}
Let  $h_{k}$ and $g_l$ be the eigenvectors associated to $\lambda_{k}(S^D_{\overset{\circ}{A}})$ and
$\lambda_l(S^D_{\overset{\circ}{B}})$ respectively, for $1 \leq k \leq \sharp V_{\overset{\circ}{A}}$ and $1\leq l\leq\sharp V_{\overset{\circ}{B}}$.\\
We define $f$ on $\mathcal{C}_{m}(\overset{\circ}{A}\oplus\overset{\circ}{B})$ by:
\begin{equation*}
f(x)=~\left\{
\begin{aligned}
&h_{k}(x)~~~~~if~x\in V_{\overset{\circ}{A}} \\
&g_l(x)~~~~~~~~if~x\in V_{\overset{\circ}{B}} \\
&0~~~~~~~~~~~~~if~x\in V_U.
  \end{aligned}
\right.
\end{equation*}
We have
$$(S f,f)_{m}=
(S^D_{\overset{\circ}{A}} h_{k},h_{k})_{m}+(S^D_{\overset{\circ}{B}} g_l,g_l)_{m}.$$

Using the Poincar\'e min-max principle (\ref{poin}), we obtain:
\begin{align*}
\lambda_{k+l}(S_G)\leq\max_{f\in F_{\overset{\circ}{A}}+F_{\overset{\circ}{B}}}\dfrac{(S_G f,f)_{m}}{( f,f)_{m}}& =\max_{h\in F_{\overset{\circ}{A}},g\in F_{\overset{\circ}{B}}}\dfrac{ (S^D_{\overset{\circ}{A}} h,h)_{m}+(S^D_{\overset{\circ}{B}} g,g)_{m}}{( h,h)_{m}+(g,g)_m}\\
&\leq \dfrac{ \lambda_{k}(S^D_{\overset{\circ}{A}})(h,h)_{m}+\lambda_l(S^D_{\overset{\circ}{B}})( g,g)_{m}}{( h,h)_{m}+(g,g)_m}\\
&\leq \max\big(\lambda_{k}(S^D_{\overset{\circ}{A}}),\lambda _l(S^D_{\overset{\circ}{B}})\big).
\end{align*}
Hence
$$\lambda_{k+l}(S_G)\leq \max\big(\lambda_{k}(S^D_{\overset{\circ}{A}}),\lambda _l(S^D_{\overset{\circ}{B}})\big).$$
\end{demo}
\begin{rem}

The previous Proposition remains true in the particular case of the Laplacian $\tilde{S}_G$.
\end{rem}
An estimate of $\lambda_{2}(S_G)$ can also be obtained with respect to $\lambda_{1}(\Delta^D_{\overset{\circ}{B}})$
and $\lambda_{1}(\Delta^D_{\overset{\circ}{B}})$.
\begin{coro}
Under the same hypothesis of the previous Proposition we have
  $$\lambda_{2}(S_G)\leq \max \Big(2\mathcal{R}e\big(\lambda_{1}(\Delta^D_{\overset{\circ}{A}})\big),2\mathcal{R}e\big(\lambda_1(\Delta^D_{\overset{\circ}{B}})\big)\Big).$$

\end{coro}
\begin{coro}
Under the same hypothesis of the previous Proposition, from the cycle graph $C_n$, we have
  $$\mathcal{R}e\big(\lambda_{2}(\Delta_{C_n})\leq \max \Big(\mathcal{R}e\big(\lambda_{1}(\Delta^D_{\overset{\circ}{A}})\big),\mathcal{R}e\big(\lambda_1(\Delta^D_{\overset{\circ}{B}})\big)\Big).$$
\end{coro}
\textbf{Acknowledgement}: I take  this  opportunity to  express  my  gratitude  to
 my PhD advisors  Colette Ann\'e and Nabila Torki-Hamza for all the fruitful discussions,
 helpful suggestions and their guidance during this work. This work was financially supported by
 the "PHC Utique" program of the French Ministry of Foreign Affairs and Ministry of higher education
 and research and the Tunisian Ministry of higher education and scientific research in the CMCU project
  number 13G1501 " Graphes, G\'eom\'etrie et th\'eorie Spectrale". Also I would like to thank the Laboratory
  of Mathematics Jean Leray of Nantes (LMJL) and the research unity (UR/13ES47) of Faculty of Sciences of
   Bizerte (University of Carthage) for its financial and its continuous support.

\end{document}